\documentclass{amsart}
\usepackage{amssymb}
\usepackage{amsmath}

\newtheorem{theorem}{Theorem}[section]
\newtheorem{lemma}[theorem]{Lemma}
\newtheorem{prop}[theorem]{Proposition}

\newtheorem{Main theorem}{Main theorem}

\theoremstyle{definition}
\newtheorem{definition}[theorem]{Definition}

\theoremstyle{remark}
\newtheorem{remark}[theorem]{Remark}

\newtheorem{acknowledgements}[theorem]{Acknowledgements}

\numberwithin{equation}{section}

\newcommand{\ZZ}{\ensuremath{\mathbb{Z}}}
\newcommand{\Z}{\ensuremath{\mathbb{Z}}}
\newcommand{\Q}{\ensuremath{\mathbb{Q}}}
\newcommand{\R}{\ensuremath{\mathbb{R}}}
\newcommand{\C}{\ensuremath{\mathbb{C}}}

\newcommand{\sur}[2]{\genfrac{}{}{0pt}{1}{#1}{#2}}
\newcommand{\surtrois}[3]{\genfrac{}{}{0pt}{1}{#1}{\genfrac{}{}{0pt}{1}{#2}{#3}}}

\newcommand{\BbbA}{\ensuremath{\mathbb{A}}}
\def\calO{{\mathcal O}}
\def\cO{{\mathcal O}}
\def\cP{{\mathcal P}}

\newcommand{\calP}{{\mathcal P}}
\newcommand{\abs}[1]{\left|{#1}\right|}
\def\grad{{\rm grad}}

\begin{document}
 \newcommand{\Xk}{\mathrm{Sym}^{k-2}(\R^2)}
 \newcommand{\Zk}{\Z_{k-2}[X,Y]}

 \title[Fundamental Theorem modulo $p^m$]{The Fundamental Theorem of prehomogeneous vector spaces
 modulo $p^m$\\ (With an appendix by F.~Sato)}

\author[R. Cluckers]{Raf Cluckers$^*$}
\address{Katholieke Universiteit Leuven, Departement wiskunde,
Celestijnenlaan 200B, B-3001 Leu\-ven, Bel\-gium. Current address:
\'Ecole Normale Sup\'erieure, D\'epartement de
ma\-th\'e\-ma\-ti\-ques et applications, 45 rue d'Ulm, 75230 Paris
Cedex 05, France} \email{cluckers@ens.fr}
\urladdr{www.dma.ens.fr/$\sim$cluckers/}
\thanks{$^*$The authors are postdoctoral fellows of the Fund for Scientific Research --
 Flanders (Belgium)(F.W.O.)}

  \author[A. Herremans]{Adriaan Herremans$^*$}
\address{Katholieke Universiteit Leuven, Departement wiskunde, Celestijnenlaan 200B, B-3001 Leuven,
 Belgium}
\email{adriaan.herremans@wis.kuleuven.ac.be}



\begin{abstract}
For a number field $K$ with ring of integers ${\mathcal O}_K$, we
prove an analogue over finite rings of the form ${\mathcal
O}_K/{\mathcal P}^m$ of the Fundamental Theorem on the Fourier
transform of a relative invariant of prehomogeneous vector spaces,
where ${\mathcal P}$ is a big enough prime ideal of ${\mathcal
O}_K$ and $m>1$. In the appendix, F.\ Sato gives an application of
the Theorems \ref{theorem A},  \ref{theorem B} and the Theorems A,
B, C in J.\ Denef and A.\ Gyoja [\emph{Character sums associated
to prehomogeneous vector spaces}, Compos.\ Math., \textbf{113}
(1998) 237--346] to the functional equation of $L$-functions of
Dirichlet type associated with prehomogeneous vector spaces.
\end{abstract}
\maketitle

 \section{Introduction}
We prove an analogue over finite rings of the Fundamental Theorem
on the Fourier transform of a relative invariant of prehomogeneous
vector spaces. In general, this Fundamental Theorem expresses the
Fourier transform of $\chi(f)$, with $\chi$ a multiplicative
(quasi-)character and $f$ a relative invariant, in terms of
$\chi(f^\vee)^{-1}$, with $f^\vee$ the dual relative invariant.
Roughly speaking, M.~Sato~\cite{SatoM} proved the Fundamental
Theorem over archimedian local fields, J.~Igusa \cite{IgusapreHVS}
over $p$-adic number fields, and  J.~Denef and A.~Gyoja
\cite{denef-gyoja} over finite fields of big enough
characteristic. In \cite{KazhPol}, the regular finite field case
is reproved. When the prehomogeneous vector space is regular and
defined over a number field $K$ we prove an analogue of the
Fundamental Theorem over rings of the form ${\mathcal
O}_K/{\mathcal P}^m,$ where ${\mathcal P}$ is a big enough prime
ideal of the ring of integers ${\mathcal O}_K$ of $K$ and $m>1$,
see Theorem \ref{theorem A}. This result is derived from the
results of \cite{denef-gyoja} by using explicit calculations of
exponential sums over the rings ${\mathcal O}_K/{\mathcal P}^m$.
\par
In \cite{Sato}, F.~Sato  introduces $L$-functions of Dirichlet
type associated to regular prehomogeneous vector spaces. In the
appendix by F.\ Sato to this paper, our results are used to obtain
functional equations for these $L$-functions
and, under extra conditions, their entireness.\\

To state the main results, we fix our notation on prehomogeneous
vector spaces. Let $(G, \rho, V)$ be a reductive prehomogeneous
vector space, meaning that $G$ is a connected complex linear
reductive algebraic group, $\rho : G \rightarrow GL(V)$ is a
finite dimensional rational representation, and $V$ has an open
$G$-orbit which is denoted by $\Omega$. Assume that $(G, \rho, V)$
has a relative invariant $0\not= f \in \C[V]$ with character $\phi
\in \mathrm{Hom}(G, \C^\times),$ that is, $f(gv) = \phi(g)f(v)$
for all $g\in G$ and $v\in V.$ We assume that $f$ is a {\it
regular} relative invariant, namely, $\Omega=V\setminus f^{-1}(0)$
is a single $G$-orbit.
 Writing $\rho^\vee : G\rightarrow GL(V^\vee)$ for the dual of
$\rho$, $(G, \rho^\vee, V^\vee)$ is also a prehomogeneous vector
space, with an open $G$-orbit which is denoted by $\Omega^\vee$,
and there exists a relative invariant $0\not= f^\vee \in
\C[V^\vee]$ whose character is $\phi^{-1}$. Then
$\Omega^\vee=V^\vee\setminus f^{\vee -1}(0)$. The map $F:= \grad
\log f$ is an isomorphism between $\Omega$ and $\Omega^\vee$ with
inverse $F^\vee:= \grad \log f^\vee$.
 One has $\dim V = \dim V^\vee =: n$  and $\deg f = \deg f^\vee =:d$.
\par

Let $K$ be a number field with ring of integers ${\mathcal O}_K$.
Suppose that $(G,\rho,V)$ is defined over $K$. We fix a basis of
the $K$-vector space $V(K)$ and we suppose that $f$ is in $K[V]$
and has coefficients in $\cO_K$ (with respect to the fixed
$K$-basis of $V(K)$). Similarly we suppose that $f^\vee$ is in
$K[V^\vee]$ and has coefficients in $\cO_K$ (with respect to the
basis of $V^\vee$ dual to the fixed basis of $V$). Write
$V(\cO_K)$ for the points of $V(K)$ with coefficients in $\cO_K$
(with respect to the fixed $K$-basis of $V(K)$), and similarly for
$V^\vee(\cO_K)$. For $I$ an ideal of $\cO$, write $V(\cO_K/I)$ for
the reduction modulo $I$ of the lattice $V(\cO_K)$.

\par
The Bernstein-Sato polynomial $b(s)$ of $f$ is defined by
$f^\vee(\grad_x)f(x)^{s+1}=b(s)f(x)^s$. Write $b_0$ for the
coefficient of the term of highest degree of $b(s)$; one has $b_0\in K$.\\

\par
The following theorem is an analogue of the Fundamental Theorem
for prehomogeneous vector spaces.
\begin{theorem}\label{theorem A}
Let $m\geq 2$ be an integer, ${\mathcal P}$ be a prime ideal of
$\cO_K$ above a big enough prime $p\in\ZZ$, $\chi$ be a primitive
multiplicative character modulo ${\mathcal P}^m$ (extended by zero
outside the multiplicative units), and let $\psi$ be a primitive
additive character modulo ${\mathcal P}^m$. Write $q:=\#
({\mathcal O}_K/{\mathcal P})$.
 For $L\in V^\vee({\mathcal O}_K/{\mathcal P}^m)$ write
\[
S(L) := \sum_{x\in V({\mathcal O}_K/{\mathcal P}^m)} \chi(f(x))
\psi(L(x)).
\]
Then the following hold:
\begin{enumerate}
\item if $f^\vee(L)\not\equiv 0 \bmod {\mathcal P}$, then
 \[
 S(L) =
q^{mn/2}\left(\frac{\sum\limits_{y\in {\mathcal O}_K/{\mathcal
P}^m} \chi^d(y) \psi(y)}{q^{m/2}}\right)
\chi(\frac{b_0f^\vee(L)^{-1}}{d^d})\alpha(\chi,m)^{n-1}\kappa^\vee(L),
 \]
where $\kappa^\vee(L)$ and $\alpha(\chi,m)$ are $1$ or $-1$;
 \item if
$f^\vee(L)\equiv 0 \bmod {\mathcal P}$, then $S(L) = 0$.
\end{enumerate}
\end{theorem}

The essential (and typical) content of this Fundamental Theorem is
that the discrete Fourier transform of the function $\chi(f)$ on
$V({\mathcal O}_K/{\mathcal P}^m)$ is equal to the function
$\chi(f^\vee)^{-1}$ on $V^\vee({\mathcal O}_K/{\mathcal P}^m)$
times some factors, and vice versa.
\par
The first part of Theorem \ref{theorem A} is obtained by combining
explicit calculations of character sums of quadratic functions (\S
2) and of discrete Fourier transforms (\S 4), a $p$-adic version
of the Lemma of Morse (\S 3), and results of \cite{denef-gyoja}.
The second part of Theorem \ref{theorem A} is established by
comparing the $L_2$-norms of $\chi(f)$ and of its discrete Fourier
transform.
\par
We also obtain explicit formulas for the constants
$\kappa^\vee(L)$ and $\alpha(\chi,m)$ of Theorem \ref{theorem A},
by using work in \cite{denef-gyoja} and elementary calculations.
To state these formulas we use the notion of the discriminant of a
matrix, as in \cite[9.1.0]{denef-gyoja}.

\begin{definition}\label{defdiscr}
For a symmetric $(n,n)$-matrix $A$ with entries in a field $k$, if
$^tXAX=$ diag $(a_1,\ldots,a_m,0,\ldots,0)$ with $X\in$GL$_n(k)$,
$^tX$ its transposed, and $a_i\in k^\times$, put
$\Delta(A):=\prod_{i=1}^m a_i\in k^\times/k^\times{}^2$, with
$k^\times{}^2$ the squares in $k^\times$, and call it the
discriminant of $A$.
\end{definition}

Write $k_{\mathcal P}$ for the finite field ${\mathcal
O}_K/{\mathcal P}$ and $k_{\mathcal P}^\times{}^2$ for the squares
in $k_{\mathcal P}^\times$. For $m>1$ and $L$ in $V^\vee({\mathcal
O}_K/{\mathcal P}^m)$ with $f^\vee(L)\not\equiv 0 \bmod {\mathcal
P}$, denote by $h^\vee(L)$ the image in $k_{\mathcal
P}^\times/k_{\mathcal P}^\times{}^2$ of the discriminant of the
matrix $\big( \frac{\partial^2 \log f^\vee} {\partial y_i\partial
y_j}(L)\big)_{ij}$, where $\{y_1,\ldots,y_n\}$ is the previously
fixed $K$-basis of $V^\vee(K)$. Write $\chi_{\frac{1}{2}}$ for the
Legendre symbol $\bmod~{\mathcal P}$. We then obtain
\begin{theorem} \label{theorem B}
The following hold in case (1) of Theorem \ref{theorem A}
\begin{itemize}
\item[(1)] $\kappa^\vee(L)=\chi_{\frac{1}{2}}(-d\, 2^{n-1}\, h^\vee(L))^m$;\\
\item[(2)] $ \alpha(\chi, m) = 1$ for $m$  even;\\
\item[(3)] $\alpha(\chi,m) = G(\chi_{\frac{1}{2}}, \psi')/\sqrt{q}$
for $m$ odd, with $\psi'$ any additive character defined by $y
\mapsto \chi(1+\pi_{\mathcal P}^{m-1}y),$ $\pi_{\mathcal P}$ any
element in ${\mathcal P}$ of ${\mathcal P}$-adic order $1$,
$\chi_{\frac{1}{2}}$ the Legendre symbol $\bmod~{\mathcal P}$, and
$G(.,.)$ the classical Gauss sum.
\end{itemize}
\end{theorem}
\begin{remark}
It is interesting to compare the formulas of Theorems \ref{theorem
A} and \ref{theorem B} to the formulas for $m=1$ given in
\cite{denef-gyoja}; it seems that for $m=1$ the formulas depend
more on subtle information of the Bernstein-Sato polynomial of
$f$.
\end{remark}
\begin{acknowledgements}
The authors would like to thank Jan Denef and Akihiko Gyoja for
the discussions and helpful ideas they communicated during the
preparation of this paper. Special thanks to Denef for giving us
an outline of the main results and their proofs. Furthermore we
are in great debt to Fumihiro Sato, who agreed on writing an
appendix to this paper. We thank the referee for his useful
suggestions to improve the presentation of the paper.
\end{acknowledgements}

\section{Preliminaries on character sums}\label{section charactersums}
Let $\chi$ be a multiplicative character $\bmod~{\mathcal P}^m,$
extended by zero for $a\equiv 0\bmod {\mathcal P}$. Say that
$\chi$ is induced by a character $\chi_1~\bmod {\mathcal P}^n$ for
$n<m$ if $\chi_1(a \bmod {\mathcal P}^n)=\chi(a)$ for all $a \in
{\mathcal O}_K/{\mathcal P}^m.$ Call $\chi$ {\it primitive}
$\bmod~{\mathcal P}^m$ if there exists no such $n<m$ and
$\chi_1~\bmod {\mathcal P}^n$ such that $\chi$ is induced by
$\chi_1$.
Analogously, call an additive character $\psi$ $\bmod~{\mathcal
P}^m$ {\it primitive} if it is not induced by a character
$\psi_1\bmod {\mathcal P}^n$ for $n<m$.
Let $f$ be a polynomial in $n$ variables over ${\mathcal O}_K.$ If
we evaluate $\sum \chi(f(x) \bmod {\mathcal P}^m,$ with $\chi$
primitive $\bmod~{\mathcal P}^m$ for some $m>1$, it is well known
that only the critical points $\bmod~{\mathcal P}$ contribute to
the sum, i.e.\ the elements $c\in {\mathcal O}_K/{\mathcal P}^m$
for which $\grad(f)|_c \equiv 0 \bmod {\mathcal P}.$ The following
is an extension of this result.
\begin{prop} \label{denefs lemma}
Let ${\mathcal P}$ be a prime ideal of ${\mathcal O}_K$, $m>1$ an
integer, $\chi$ a primitive multiplicative character
$\bmod~{\mathcal P}^m$ and $f\in{\mathcal O}_K[X]$ a polynomial in
$n$ variables. Writing
 $$
 S_f = \sum_{x\in
({\mathcal O}_K/{\mathcal P}^m)^n}\chi(f(x)), $$ we obtain that
$$ S_f = \sum_{\sur{x\in ({\mathcal O}_K/{\mathcal P}^m)^n}
{v_{{\mathcal P}}(\grad(f)|_x)\geq {(m-1)/ 2}}} \chi(f(x)), $$
where $v_{{\mathcal P}}(\grad(f)|_x)$ is the minimum of the
${\mathcal P}$-valuations of $\frac{\partial f}{\partial x_j}|_x$
for $j=1,\ldots,n$. Moreover, the same formulas hold for an
additive character $\psi$ instead of $\chi$.
\end{prop}
\begin{proof} We treat the case that $f$ is a function of one
variable; the general case is completely analogous. It suffices to
prove that
$$
S_f(i,c) := \sum_{\surtrois{x\in~{\mathcal O}_K/{\mathcal P}^m }{
v_{{\mathcal P}}(f'(x)) = i }{ x\equiv c\bmod {\mathcal
P}^{m-i-1}}} \chi(f(x))$$ is zero for every $i<{(m-1)/ 2}$ and
every $c \in {\mathcal O}_K/{\mathcal P}^{m-i-1}.$ If $f(c) \equiv
0 \bmod {\mathcal P},$ then $S_f(i,c)$ is trivially zero. Suppose
$f(c)\not\equiv 0 \bmod {\mathcal P}.$ Let $\pi_{\mathcal P}$ be
in ${\mathcal P}$ of ${\mathcal P}$-adic order $1$. Writing $f(x)$
as a polynomial in $(x-c)$, we get the following equalities
\begin{eqnarray}
S_f(i,c) &=& \sum\limits_{\surtrois{x\in~{\mathcal O}_K/{\mathcal
P}^m} { v_{{\mathcal P}}(f'(x)) = i }{
  x\equiv c\bmod {\mathcal P}^{m-i-1}}} \chi(f(c) + (x-c)
f'(c) + \ldots) \label{1a}
\\
&=&
 \chi(f(c))
\sum\limits_{\surtrois{x\in~{\mathcal O}_K/{\mathcal P}^m
}{v_{{\mathcal P}}(f'(x)) = i }{ x\equiv c\bmod {\mathcal
P}^{m-i-1}}} \chi\left(1 + (x-c) \frac{f'(c)}{ f(c)}\right)
\label{1b}
\\
  &=& q^i\chi(f(c)) \sum_{z\in~{\mathcal O}_K/{\mathcal P}}
\chi(1+\pi_{\mathcal P}^{m-1} z\alpha)\label{1d}\\ &=& 0.
\nonumber
\end{eqnarray} Equality (\ref{1b}) comes from the fact that the
${\mathcal P}$-valuation of the terms of degree $\geq 2$ in $x-c$
are at least $m.$ Indeed $2m-2i-2\geq m$ if $i\leq {(m-2)/ 2}$ or
$i<{(m-1)/ 2}.$
 Either there are no terms in (\ref{1b}) and
then it is automatically zero, or, there are terms and then
Equality (\ref{1d}) follows immediately. We conclude that the sum
(\ref{1d}) is zero since we sum a non-trivial additive character
$\psi(z):=\chi(1+\pi_{\mathcal P}^{m-1}\alpha z)$ over ${\mathcal
O}_K/{\mathcal P}$; $\psi$ is indeed non-trivial since $\chi$ is
primitive.
\par
\end{proof}

\begin{definition} \label{definitie alpha(chi,m)}Let
${\mathcal P}$ be a prime ideal of ${\mathcal O}_K$ not containing
2 and let $m>1$ be an integer. Let $\chi$, resp.\ $\psi$, be a
primitive multiplicative, resp.\ primitive additive, character
$\bmod\ {\mathcal P}^m.$ Then, put $$ \tilde\alpha(\chi, m) :=
\sum_{x\in~{\mathcal O}_K/{\mathcal P}^m} \chi(1 + x^2),
$$ and
$$ \tilde\alpha(\psi, m) :=
\sum_{x\in~{\mathcal O}_K/{\mathcal P}^m} \psi(x^2). $$
\end{definition}
One can calculate the value of $\tilde\alpha(\chi, m)$, using
Proposition \ref{denefs lemma} in an elementary way, to obtain the
following lemma.
\begin{lemma} \label{berekening alpha(chi,m)}
With the assumptions and notation of Definition \ref{definitie
alpha(chi,m)}, the following hold
\begin{itemize}
 \item $ \tilde\alpha(\chi, m) = q^{{m/
2}}$ if $m$ is even;
 \item $\tilde\alpha(\chi,m) = q^{{(m-1)/
2}} G(\chi_{\frac{1}{2}} , \psi')$ if $m$ is odd. Herein,
$\chi_{\frac{1}{2}}$ is the Legendre symbol $\bmod~{\mathcal P}$,
$\psi'$ is any additive character defined by $y \mapsto
\chi(1+\pi_{\mathcal P}^{m-1}y)$ with $\pi_{\mathcal P}$ any
element in ${\mathcal P}$ of ${\mathcal P}$-adic order $1$, and
$G(.,.)$ is the classical Gauss sum.
\end{itemize}
\end{lemma}

This lemma and an induction argument yield the following
proposition.

\begin{prop} \label{chi van een kwadratische functie}
Use the assumptions and notation of Definition \ref{definitie
alpha(chi,m)}. Suppose that $f=a_0+a_1x^2_1+\ldots+a_nx^2_n,$
where $a_0,a_1,\ldots,a_n$ are multiplicative units
$\bmod~{\mathcal P}^m$. Putting
$$
S_f := \sum\limits_{x\in\left({\mathcal O}_K/{\mathcal
P}^m\right)^n} \chi (a_0+a_1x^2_1+\ldots+a_nx^2_n),
 $$
the following holds
 $$
 S_f = \chi(a_0)\chi_{\frac{1}{2}}
 (a^n_0a_1\ldots a_n)^m
\tilde\alpha(\chi,m)^n
 $$
  with $\chi_{\frac{1}{2}}$
   the Legendre
symbol $\bmod\ {\mathcal P}.$
\end{prop}

Similar formulas can be obtained for $ \tilde\alpha(\psi, m)$ with
$\psi$ a primitive additive character modulo ${\mathcal P}^m$.
%
\section{A $p$-adic analogue of the Lemma of Morse}
We prove an analogue of the Lemma of Morse (originally formulated
for $C^{\infty}$ functions on real manifolds), for $p$-adic
analytic functions. This Lemma is normally stated as a local
property, but here we can work with fixed (large) neighborhoods.
To deal with the fact that $\Z_p$ is totally disconnected, we will
use a global notion of analycity for $p$-adic maps. Results in
this section also hold for finite field extensions of $\Q_p$.

\begin{definition}
Let $A\subset \Z_p^n$ be open. Call a function $f:A\to\Z_p$ {\it
globally analytic} if there is a power series $\sum
c_ix^i\in\Z_p[[x]]$ which converges on $A$ such that $f(x)=\sum
c_i x^i$ for each $x\in A$, with $x=(x_1,\ldots,x_n)$,
$i=(i_1,\ldots,i_n)$ and $x^i=\prod_{j=1}^nx_j^{i_j}$. Call $f$
analytic if there is a finite cover of $A$ by opens $U$ such that
the restriction of $f$ to each $U$ is globally analytic.
Similarly, call a map $f:A\subset \Z_p^n \to \Z_p^m$ analytic if
it is given by analytic functions on $A$. An analytic bijection
with analytic inverse is called {\it bi-analytic}. For an analytic
function $f:A\subset \Z_p^n \to\Z_p$, define the gradient of $f$
in $a\in A$ as $\grad(f)|_{a}=(\frac{\partial f}{\partial
x_1}|_{a},\ldots,\frac{\partial f}{\partial x_n}|_{a})$ and the
Hessian of $f$ in $a$ as
$\mathrm{Hs}(f)|_{a}=(\frac{\partial^2f}{\partial x_i\partial
x_j}|_{a})_{ij}$. Say that $a$ is a critical point if
$\grad(f)|_{a}=0$ and call the critical point $a\in A$
non-degenerate if $\det(\mathrm{Hs}(f)|_{a})\not=0$. Call $a\in A$
a {\it non-degenerate critical point modulo $p$} if
$\grad(f)|_{a}\equiv0\bmod p$ and
$\det(\mathrm{Hs}(f)|_{a})\not\equiv0\bmod p$.\\
Let $M\subset \Z_p^n$ be a compact $p$-adic manifold of pure
dimension $d$ and $f:M\to \Z_p$ an analytic function. Then there
exists a finite disjoint cover of $M$ by opens $U_i$ and analytic
isometries $\pi_i$ from $U_i$ to open balls in $\Z_p^d$ such that
the maps $f\circ\pi^{-1}_i$ are globally analytic, see
e.g.~\cite{serre}. Call a point $a\in U_i\subset M$ a critical
point of $f$ if $\pi_i(a)$ is a critical point of
$f\circ\pi^{-1}_i$ as defined above. Similarly we speak of
non-degenerate critical points and non-degenerate critical point
modulo $p$. This is independent of the choice of $U_i$ and
$\pi_i$.
\end{definition}
The next Lemma is a $p$-adic variant of the Inverse Function
Theorem.
\begin{lemma}\cite[Cor. 2.2.1]{I}\label{inversefunctie}
Suppose that $g_1,\ldots,g_n\in\Z_p[[x_1,\ldots,x_n]]$ satisfy
$\det(\frac{\partial g_i}{\partial x_j}|_0)_{ij}\not\equiv0 \bmod
p$ and $g_i(0)\equiv 0\bmod p$ for $i=1,\ldots,n$. Then the map
$g:(p\Z_p)^n\to(p\Z_p)^n:x\mapsto(g_1(x),\ldots,g_n(x))$ is
(globally) bi-analytic.
\end{lemma}
The following is a $p$-adic analogue of the Lemma of Morse. The
proof goes along the same lines as in \cite{M} and we refer to
\cite{M}, Lemma 2.2 for the details.
\begin{lemma}[Morse]\label{Morse}
Let $p\not=2$ and let $f:(p\Z_p)^n\to p\Z_p$ be a globally
analytic map (thus $f$ is given by a single power series in
$\Z_p[[x]]$) such that $0$ is a non-degenerate critical point
modulo $p$. Then there is a unique critical point $c$ of $f$ and
there is a (globally) bi-analytic map
$T:(p\Z_p)^n\to(p\Z_p)^n:x\mapsto T(x)=u$ such that
\[f(x)=f(c)+\sum_{i=1}^n a_iu_i^2,\ \mbox{ for all }x\in(p\Z_p)^n,\]
with $a_i$ in $\Z_p^\times$. Moreover,
$\chi_{\frac{1}{2}}(H)=\chi_{\frac{1}{2}}(\prod_ia_i)$ with
$H=2^{-n}\det(\mathrm{Hs}(f)|_c)$.
\end{lemma}
\begin{proof} The uniqueness of the critical point $c$ is proved
in Lemma \ref{crit} below. We may suppose that $f(0)=0$ and $c=0$.
We can write $f(x)=\sum_{i,j}x_ix_jh_{ij}(x)$ with
$h_{ij}(x)\in\Z_p[[x]]$ and since $p\not=2$ we can assume that
$h_{ij}=h_{ji}$. Suppose by induction that we have a (globally)
bi-analytic map $T_{r-1}:(p\Z_p)^n\to (p\Z_p)^n:x\mapsto u$ such
that for each $x\in (p\Z_p)^n$
\[f(x)=\sum_{i=1}^{r-1}a_iu_i^2+\sum_{i,j\geq r}u_iu_jH_{ij}(u),\]
with $H_{ij}(u)\in\Z_p[[u]]$, $H_{ij}=H_{ji}$ and
$a_i\in\Z_p^\times$. We have
$\det(\mathrm{Hs}(f)|_0)=\lambda\det(H_{ij}(0))_{ij}$ with
$\lambda$ a unit in $\Z_p$, so $\det(H_{ij}(0))_{ij}\not\equiv
0\bmod p$ and after a linear change in the last $n-r+1$
coordinates we may assume that $H_{rr}(0)\not\equiv 0\bmod p$. Put
$a_r= H_{rr}(0).$ Then $H_{rr}(u)/a_r\equiv1\bmod p$ for each
$u\in (p\Z_p)^n$. There is a well-defined square root function
$\sqrt{}:1+p\Z_p\to 1+p\Z_p$ which is (globally) analytic. Put
$g(u)=\sqrt{\frac{H_{rr}(u)}{a_r}}$ for $u\in (p\Z_p)^n$. Then,
$g$ is (globally) analytic on $(p\Z_p)^n$. We can now define a
(globally) analytic map $T'_{r-1}:(p\Z_p)^n\to(p\Z_p)^n:u\mapsto
T'_{r-1}(u)=(v_i)$ by
\begin{eqnarray*}
v_i & = & u_i\qquad \mbox{ if }i\not=r,\\ v_r & = &
g(u)\big(u_r+\sum_{i>r}u_i\frac{H_{ir}(u)}{H_{rr}(u)}\big).
\end{eqnarray*}
Clearly we have that $\det(\frac{\partial v_i}{\partial
u_j}|_0)=g(0)\not\equiv 0\bmod p$, thus, the map $T'_{r-1}$ is
(globally) bi-analytic by Lemma \ref{inversefunctie}. Put
$T_r=T'_{r-1}\circ T_{r-1}.$ We then obtain for all
$x\in(p\Z_p)^n$ and $v=T_r(x)$
\[f(x)=\sum_{i=1}^{r}a_iv_i^2+\sum_{i,j> r}v_iv_jH'_{ij}(v),\]
with $H'_{ij}(v)\in\Z_p[[v]]$. This finishes the induction
argument.
 \par
The equality
$\chi_{\frac{1}{2}}(\det(\mathrm{Hs}(f)|_0)=\chi_{\frac{1}{2}}(\prod_i2a_i)$
follows by a classical argument as in \cite{M}. To finish the
proof we only have to prove the Lemma below. \end{proof}
\begin{lemma}\label{crit} Let $p\not=2$ and let $f:(p\Z_p)^n\to p\Z_p$ be a globally
analytic map (i.e., given by a power series in $\Z_p[[x]]$) such
that $0$ is a non-degenerate critical point modulo $p$. Then $f$
has a unique critical point $c\in (p\Z_p)^n$ and this is a
non-degenerate critical point.
\end{lemma}
\begin{proof}
We can write $f(x)=f(0)+\sum_ia_ix_i+g(x)$ with
$g(x)=\sum_{ij}x_ix_jh_{ij}(x)$ for some $h_{ij}\in\Z_p[[x]]$.
Since $0$ is a non-degenerate critical point of $f$ modulo $p$, we
see that $a_i\equiv0\bmod p$ and
$\det(\mathrm{Hs}(g)|_0)\not\equiv0\bmod p$. By Lemma
\ref{inversefunctie} the map
$T:(p\Z_p)^n\to(p\Z_p)^n:x\mapsto\grad(g)|_x=(\frac{\partial
g}{\partial x_1}(x),\ldots,\frac{\partial g}{\partial x_n}(x))$ is
a bi-analytic bijection. A fortiori, there is a unique point
$c\in(p\Z_p)^n$ such that $T(c)=(-a_1,\ldots,-a_n)$. Since
$\grad(f)|_x=\grad(g)|_x+(a_1,\ldots,a_n)$, the condition
$T(x)=(-a_1,\ldots,-a_n)$ is equivalent with the condition
$\grad(f)|_x=0$ for $x\in(p\Z_p)^n$. Thus $c$ is the unique
critical point of $f$. Moreover, $\det(\mathrm{Hs}(f)|_c)\not=0$
since
$\det(\mathrm{Hs}(f)|_c)\equiv\det(\mathrm{Hs}(f)|_0)\not\equiv0\bmod
p$.
\end{proof}

Finally, we give an application of the Lemma of Morse to calculate
character sums of polynomials.

\begin{prop}\label{eindig veel kritische} Let $p\not=2$,
let $f\in\Z_p[x_1,\ldots,x_n]$ be a polynomial, and let $X$ be a
smooth subvariety of $ {\mathbb{A}}^n_{{\Z_p}}$ (smooth variety
meaning a separated, reduced, irreducible scheme of finite type
and smooth over $\Z_p$). Let $M$ be the compact $p$-adic variety
$X(\Z_p)\setminus f^{-1}(p\Z_p)$ and write $d$ for the dimension
of  $M$. Let $\chi$ be a primitive multiplicative character
$\bmod~p^m$, $m>1$. Let $c_1,\ldots, c_l\in M$ be the critical
points of $f|_M:M\to \Z_p$ and suppose that the critical points
$c_1,\ldots, c_l$ are non-degenerate modulo $p$. Put
 \[
 S_f=\sum\limits_{x\in X(\Z_p/p^m\Z_p)}\chi(f(x)).
 \]
Then the following holds
 \[
S_f=\sum\limits_{i=1}^l
 \chi(f(c_i))\chi_{\frac{1}{2}}(f(c_i)^dH_{c_i})^m\tilde\alpha(\chi,m)^d,
  \]
with $H_{c_i}=2^{-d}\Delta(\mathrm{Hs}(f)_{|c_i})$,
$\tilde\alpha(\chi,m)$ as in section \ref{section charactersums},
and $\Delta$ the discriminant as defined in \ref{defdiscr}.
\end{prop}

\begin{proof}
By Hensels Lemma and since $X$ is smooth over $\Z_p$, we can cover
$M$ by finitely many disjoint compact opens $U$ such that for each
$U$ we can find an analytic isometry $\pi:U\to(p\Z_p)^d$. For
$i=1,\ldots,r$, let $U_{i}$ be the open in this cover containing
$c_i$ and we write $\pi_i$ for the corresponding isometry. By
Proposition \ref{denefs lemma} and Lemma \ref{crit} it follows
that $S_f=\sum\limits_{i=1}^l S_f(c_i)$ where
 \[
 S_f(c_i)=\sum\limits_{x\in(p\Z_p/p^m\Z_p)^d}
        \chi(f\circ\pi_i^{-1}(x)).
  \]
By the Lemma of Morse we can find for each critical point $c_i$ a
bi-analytic isometric transformation $T_{i}: (p\Z_p)^d
\to(p\Z_p)^d$ such that
\[f(x)=f(c_i)+\sum_{j=1}^d a_ju_j^2\ \mbox{ for all }x\in U_i \mbox{ and }
u=T_{i}(\pi_i(x)),\] with $a_j$ in $\Z_p^\times$ and
$\chi_{\frac{1}{2}} (H_{c_i})=\chi_{\frac{1}{2}} (\prod_j a_j)$.
We can now calculate
\begin{eqnarray}
 S_f(c_i)
      & = &
      \sum\limits_{u\in(p\Z/p^m\Z)^d}
            \chi(f(c_i)+\sum_{j=1}^d a_ju_j^2)\label{s3}
\\
     & = &
      \sum\limits_{u\in(\Z/p^m\Z)^d}  \chi(f(c_i)+\sum_{j=1}^d a_ju_j^2)\label{s'3}
\\
     & = &
      \chi(f(c_i))\chi_{\frac{1}{2}}(f(c_i)^d\prod_ja_j)^m\tilde\alpha(\chi,m)^d\label{s4}
\\
     & = &
      \chi(f(c_i))\chi_{\frac{1}{2}}(f(c_i)^dH_{c_i})^m\tilde\alpha(\chi,m)^d\label{s5}.
\end{eqnarray}

Equality (\ref{s3}) is clear. Also equality (\ref{s'3}) is easy
and follows by similar arguments as in the proof of Proposition
\ref{denefs lemma}. Equality (\ref{s4}) comes from Proposition
\ref{chi van een kwadratische functie} and the last equality holds
because $\chi_{\frac{1}{2}} (H_{c_i})=\chi_{\frac{1}{2}} (\prod_j
a_j)$.
\end{proof}

\section{Discrete Fourier transforms of characters of homogeneous polynomials}
Let ${\mathcal P}$ be a prime ideal of ${\mathcal O}_K$ and let
$\pi_{\mathcal P}\in{\mathcal P}$ be of ${\mathcal P}$-adic order
$1$. Let $L(x) = \sum\limits^n_{i=1} a_ix_i$ be a linear form on
$({\mathcal O}_K/{\mathcal P}^m)^n,$ with $m>1$ and $a_i\in
{\mathcal O}_K/{\mathcal P}^m$. Let $f \in {\mathcal O}_K[x]$ be a
homogeneous polynomial of degree $d$ and let $\chi$, resp.\
$\psi$, be a primitive multiplicative, resp.\ primitive additive,
character $\bmod {\mathcal P}^m.$

We will calculate the discrete Fourier transform of $\chi(f)$,
defined by
 $$
 S(L) := \sum_{x\in ({\mathcal O}_K/{\mathcal P}^m)^n} \chi(f(x))
 \psi(L(x)).
 $$

After a linear change of variables one can assume that the
${\mathcal P}$-valuation of $a_1$ is minimal among the
$v_{{\mathcal P}}(a_i)$. Write $k:=v_{{\mathcal P}}(a_1)$ and
$L(x)=\pi_{\mathcal P}^k (a_1'x_1 + a_2'x_2+\ldots+a_n'x_n)$ for
some $a_i'\in {\mathcal O}_K/{\mathcal P}^m$. After applying the
invertible linear transformation $$(x_1, \ldots, x_n) \mapsto
(a_1'x_1 + \ldots + a_n'x_n, x_2, \ldots ,x_n),$$ one reduces to
the case that $L(x)=\pi_{\mathcal P}^k x_1$. This reduction is
used in the proof of the following result.

\begin{prop} \label{fourier} Let ${\mathcal P}$ be a prime ideal of ${\mathcal O}_K$ and let
$\pi_{\mathcal P}\in{\mathcal P}$ be of ${\mathcal P}$-adic order
$1$. Let $L(x) = \sum\limits^n_{i=1} a_ix_i$ be a linear form on
$({\mathcal O}_K/{\mathcal P}^m)^n,$ with $m>1$ and $a_i\in
{\mathcal O}_K/{\mathcal P}^m$. Let $f \in {\mathcal O}_K[x]$ be a
homogeneous polynomial of degree $d$ and let $\chi$, resp.\
$\psi$, be a primitive multiplicative, resp.\ primitive additive,
character $\bmod {\mathcal P}^m.$ Put $k:=
\min\limits_{i=1,\ldots,n} v_{{\mathcal P}}(a_i)$ and put $$
 S(L)=
\sum\limits_{x\in ({\mathcal O}_K/{\mathcal P}^m)^n} \chi(f(x))
\psi(L(x)). $$ If ${\mathcal P}\nmid d$, then
 $$
 S(L) = \left\{
 \begin{array}{cl}
 0& \mbox{if $k\not=0$;}\\
 \left(\sum\limits_{y\in {\mathcal
O}_K/{\mathcal P}^m} \chi^d(y)
 \psi(y)\right) \left( \sum\limits_{\sur{x \in
 ({\mathcal
O}_K/{\mathcal P}^m)^{n} }{
 L(x)\equiv 1\bmod {\mathcal P}^m }} \chi(f(x)) \right)& \mbox{if $k=0.$}
 \end{array}\right.
 $$
\end{prop}
\begin{proof}
As explained in the previous discussion, we may assume that
$L(x)=\pi_{\mathcal P}^k x_1,$ with $\pi_{\mathcal P}\in{\mathcal
P}$ of ${\mathcal P}$-adic order $1$. We split up the sum
depending on the ${\mathcal P}$-valuation of $x_1.$ We call the
subsum of $S(L),$ consisting only of the elements $x$ with
$v_{\mathcal P}(x_1)=j$ to be $A_{jk}.$ More precisely, we have
$S(L)=\sum\limits^{m}_{j=0} A_{jk},$ with
$$A_{jk} := \sum_{\surtrois{y\in ({\mathcal O}_K/{\mathcal P}^{m-j})^\times
}{ x_1 = y\pi_{\mathcal P}^j}{\widehat{x}\in ({\mathcal
O}_K/{\mathcal P}^m)^{n-1}}} \chi(f(y\pi_{\mathcal P}^j,
\widehat{x}))
 \psi(\pi_{\mathcal P}^{j+k}y),$$
where $\widehat{x}:=(x_2,\ldots,x_n).$ Rewrite $q^jA_{jk}$ as
\begin{eqnarray}
&& \sum_{\sur{y\in ({\mathcal O}_K/{\mathcal P}^m)^\times}
{\widehat{x}\in ({\mathcal O}_K/{\mathcal P}^m)^{n-1}}}
\chi(f(y\pi_{\mathcal P}^j, \widehat{x})) \psi(\pi_{\mathcal
P}^{j+k}y) \label{3a}
\\
&=&
 \sum_{\sur{y\in ({\mathcal O}_K/{\mathcal P}^m)^\times}{
\widehat{x}\in ({\mathcal O}_K/{\mathcal P}^m)^{n-1} }}
\chi(f(y\pi_{\mathcal P}^j, y\hat x))
  \psi(\pi_{\mathcal P}^{j+k}y) \label{3b}
  \\
&=&
 \left(\sum\limits_{y\in {\mathcal O}_K/{\mathcal P}^m} \chi^d(y) \psi(\pi_{\mathcal P}^{j+k}y)
 \right) \left( \sum\limits_{\widehat{x} \in ({\mathcal
O}_K/{\mathcal P}^m)^{n-1}} \chi(f(\pi_{\mathcal
P}^j,\widehat{x})) \right) \label{3c}
\end{eqnarray}
Equality (\ref{3a}) holds because only the value of $y \bmod
{\mathcal P}^{m-j}$ is relevant. Equality (\ref{3b}) is just
substituting $\widehat{x} = (x_2, \ldots, x_n)$ by $y\widehat{x} =
(yx_2, \ldots, yx_n)$; since $y$ is a unit the set over which we
sum does not change.
 The last equality uses that $f$ is homogeneous of degree $d$ and the fact that
$\chi^d(y)=0$ if $y\not\in ({\mathcal O}_K/{\mathcal
P}^m)^\times.$
\par
 We want to prove that all $A_{jk}$ are zero except when $k=j=0.$
Since ${\mathcal P}\nmid d$, we have that $\chi^d$ is still a
primitive character $\bmod {\mathcal P}^m.$ Therefore, there
exists an $a\in {\mathcal O}_K/{\mathcal P}^m$ such that $a\equiv
1\bmod {\mathcal P}^{m-1}$ and $\chi^d(a)\not=1$ (see section
\ref{section charactersums}). By a classical argument we obtain:
\begin{eqnarray*}
\sum_{y\in {\mathcal O}_K/{\mathcal P}^m} \chi^d(y)
\psi(\pi_{\mathcal P}^{j+k}y)
  &=&
\sum_{y\in {\mathcal O}_K/{\mathcal P}^m} \chi^d(ay) \psi(a\pi_{\mathcal P}^{j+k}y)\\
 &=&
 \chi^d(a) \sum_{y\in {\mathcal O}_K/{\mathcal P}^m} \chi^d(y)
\psi(\pi_{\mathcal P}^{j+k}y),
\end{eqnarray*}
if $j+k\geq 1.$ Indeed, the first equation is just substituting
$y$ by $ay$, where $a$ is a unit. The second uses the fact that
$\chi^d$ is multiplicative and since $a\equiv 1 \bmod {\mathcal
P}^{m-1},$ we have $ay\equiv y \bmod {\mathcal P}^{m-k-j}$ if
$j+k\geq 1.$\\ Since $\chi^d(a)\not=1,$ we conclude that
 $$
 \sum_{y\in {\mathcal O}_K/{\mathcal P}^m} \chi^d(y) \psi(\pi_{\mathcal P}^{j+k}y) = 0.
 $$
Note that when $j+k\geq m$ this sum is just $\sum\limits_{y\in
{\mathcal O}_K/{\mathcal P}^m} \chi^d(y)$ which is directly seen
to be zero. This proves the Proposition since only $A_{00}$ is
non-zero. $A_{00}$ equals the non-zero term of the Proposition by
equality (\ref{3c}). \end{proof}

\section{Applications to Prehomogeneous Vector Spaces.}\label{sec:applic}

We use the notation and the asumptions from the introduction.

 \begin{proof}[Proof of Theorem \ref{theorem A}]
For any $v^\vee\in V^\vee$ let $H(v^\vee)$ be the hyperplane in
$V$ defined by $v^\vee(x)-1$. By \cite[Lemma 9.1.2]{denef-gyoja}
and \cite[Lemma 9.1.7]{denef-gyoja}, for any $v^\vee$ with
$f^\vee(v^\vee)\not=0$, the restriction $f_{|H(v^\vee)\cap\Omega}$
has the point $d^{-1}F^\vee(v^\vee)$ as its unique critical point
and it is a non-degenerate critical point. The fact that this
holds for all $v^\vee$ with $f^\vee(v^\vee)\not=0$ is a first
order statement (in the language of rings), and hence, the
analogue statement is true over finite fields of big enough
characteristic. That is, for all prime ideals $\cP$ of $\cO_K$
above big enough primes $p\in\ZZ$ with residue field $k_\cP$, for
all $v^\vee\in V^\vee(k_\cP)$ with $f^\vee(v^\vee)\not=0$, the
point $d^{-1}F^\vee(v^\vee)$ is the unique critical point of
$f_{|H(v^\vee)\cap\Omega(k_\cP)}$ and it is a non-degenerate
critical point.

Take a prime ideal $\cP$ above a big enough prime $p$. Choose
$L\in V^\vee({\mathcal O}_K/{\mathcal P}^m)$ with
$f^\vee(L)\not\equiv0 \bmod {\mathcal P}$. Let $L_0\in
V^\vee({\mathcal O}_K)$ lie above $L$. Let $c$ be
$d^{-1}F^\vee(L_0)$. Let $R$ be the valuation ring of the
$\cP$-adic completion of $K$. Since $\cP$ is supposed to lie above
a big enough prime $p\in\ZZ$, we may suppose that $d$ is
invertible in $R$, and hence, $c$ lies in $V(R)$. Make the plane
$H(L_0)$ into a vector space by choosing a zero point in
$H(L_0)(\cO_K)$. Take a basis of the lattice $H(L_0)(\cO_K)$ in
$H(L_0)$ and take coordinates on $H(L_0)$ with respect to this
basis. Let $H_{c}$ be $2^{-n+1}$ times the determinant of the
Hessian in $c$ of $f_{|H(L_0)}$, expressed in these coordinates.
Since $p\in\ZZ$ is big enough and by the above discussion, $H_c$
lies in $R$ and $H_{c}\bmod \cP$ is nonzero in $\mathcal
{O}_K/{\mathcal P}$, and this is so uniformly in $L\in
V^\vee({\mathcal O}_K/{\mathcal P}^m)$.

Then,
\begin{eqnarray}
S(L) & = &\left(\sum\limits_{y\in {\mathcal O}_K/{\mathcal P}^m}
\chi^d(y) \psi(y)\right) \left( \sum\limits_{\sur{x \in
V({\mathcal O}_K/{\mathcal P}^m) }{ L(x)\equiv 1\bmod {\mathcal
P}^m }} \chi(f(x)) \right)\label{preh1}
 \\
 & = &
\left(\sum\limits_{y\in {\mathcal O}_K/{\mathcal P}^m} \chi^d(y)
\psi(y)
  \right)\chi(f(c))\chi_{\frac{1}{2}}(f(c)^{n-1}H_{c})^m\tilde\alpha(\chi,m)^{n-1}\label{preh2}\\
 & = &
 q^{mn/2}\left(\frac{\sum\limits_{{\mathcal
O}_K/{\mathcal P}^m} \chi^d(y) \psi(y)}{q^{m/2}}\right)
\chi(\frac{b_0f^\vee(L)^{-1}}{d^d}) \kappa^\vee(L)
\frac{\tilde\alpha(\chi,m)^{n-1}}{q^{(n-1)m/2}},\label{preh2'}
\end{eqnarray}
with $\kappa^\vee(L)=\chi_{\frac{1}{2}}(f(c)^{n-1}H_{c})^m$.
Equality (\ref{preh1}) follows from Proposition \ref{fourier}.
 We obtain (\ref{preh2}) from Proposition \ref{eindig veel
kritische}. If we now use the fact that $f(c)=d^{-d} f(F^\vee(L))
=d^{-d}b_0 f^\vee(L)^{-1}$ (see \cite[Lemma 9.1.2]{denef-gyoja}
and \cite[Lemma 1.8]{gyoja}), equation (\ref{preh2'}) and thus
case (1) of the Theorem follows with
$\kappa^\vee(L)=\chi_{\frac{1}{2}}(f(c)^{n-1}H_{c})^m$ and
$\alpha(\chi,m)=\frac{\tilde\alpha(\chi,m)^{n-1}}{q^{(n-1)m/2}}$.
\par

Now we prove case (2) of Theorem \ref{theorem A} with a technique
of N.~Kawanaka. The function $S:\mathrm{Hom}(({\mathcal
O}_K/{\mathcal P}^m)^n,{\mathcal O}_K/{\mathcal P}^m)\to \C:
L\mapsto S(L)$ is the discrete Fourier transform of
$\chi(f):({\mathcal O}_K/{\mathcal P}^m)^n\to\C:x\mapsto
\chi(f(x))$. By a classical result on $L_2$-norms of Fourier
transforms on finite abelian groups, it follows that
\begin{equation}
 \|S\|_2^2 = q^{mn} \|\chi(f)\|_2^2 = q^{mn} N_1,
\label{L2}
 \end{equation}
with $N_1:=\#A$ and $A:=\{x\in V({\mathcal O}_K/{\mathcal
P}^m)\mid f(x)\not\equiv0\bmod {\mathcal P}\}.$

It follows from the formula in case (1), that $|S(L)|^2=q^{mn}$
for $L\in A$. Writing $B:=V({\mathcal O}_K/{\mathcal
P}^m)\setminus A$, one has
\begin{equation}\label{eq22}
\|S\|_2^2  =  \sum\limits_{A} q^{mn}+ \sum\limits_{B}|S(L)|^2,
\end{equation}
and hence, $S(L)=0$ for $L\in B$. \end{proof}

\begin{proof}[Proof of Theorem \ref{theorem B}]
The statement about $\alpha(\chi,m)$ follows from equation
(\ref{preh2'}) and Lemma \ref{berekening alpha(chi,m)}. We recall
that $h^\vee(L)$ is defined in the introduction, and that the
number $\chi_{\frac{1}{2}}(h^\vee(L))$ is well-defined. We obtain
the value of $\kappa^\vee(L)$ as an immediate corollary of
\cite[Lemma 9.1.7]{denef-gyoja}.
\end{proof}

\section*{APPENDIX: $L$-functions of prehomogeneous vector spaces
       (by Fumihiro Sato)}
\setcounter{section}{6} \setcounter{equation}{0} In this note, we
give an application of Theorems \ref{theorem A} and \ref{theorem
B} in \cite{CH} (and Theorems A, B, C in \cite{denef-gyoja}) to
the functional equation of $L$-functions of Dirichlet type
associated with prehomogeneous vector spaces, which is a
generalization of Theorem L in \cite{Sato}.

In the following we retain the notation in \cite{CH}. However, for
simplicity, we assume that $K=\Q$, $\calO_K=\Z$ and $\calP=(p)$
with a rational prime $p$. Let $\chi$ be a primitive Dirichlet
character with conductor $N>1$. We extend $\chi$ to $\Z/N\Z$ by
$\chi(a)=0$ for $a \not\in (\Z/N\Z)^\times$. Put $m(p)=v_p(N)$,
the $p$-order of $N$, for any rational prime $p$. Since
$(\Z/N\Z)^\times$ is isomorphic to $\prod_{p|N}
(\Z/p^{m(p)}\Z)^\times$, $\chi$ induces a primitive character
$\chi^{(p)}: (\Z/p^{m(p)}\Z)^\times \rightarrow \C^\times$ for
each $p|N$.

Let $(G,\rho,V)$ be a reductive prehomogeneous vector space
defined over $\Q$. Let $P_1,\ldots,P_\ell$ be the fundamental
relative invariants over $\Q$, namely, the $\Q$-irreducible
relatively invariant polynomials on $V$. We denote by $\phi_i$ $(1
\leq i \leq \ell)$ the rational character of $G$ corresponding to
$P_i$. The fundamental relative invariants are determined uniquely
up to a non-zero constant multiple in $\Q^\times$ and any relative
invariant in $\Q[V]$ is a monomial of them.

We fix a basis of the $\Q$-vector space $V(\Q)$ and take a
relative invariant  $f \in \Q[V]$ with coefficients in $\Z$ (with
respect to the fixed $\Q$-basis of $V(\Q)$). The character $\phi$
corresponding to $f$ is defined over $\Q$.

In the following we assume that
\begin{enumerate}
\def\labelenumi{(A.\theenumi)}
 \item $f$ is a {\it regular} relative invariant, namely,
$\Omega=V\setminus f^{-1}(0)$ is a single $G$-orbit; \item for
every $x \in \Omega(\Q)$, the group of $\Q$-rational characters of
the identity component of
\[
G_x=\left\{g \in G \left| \rho(g)x=x \right.\right\}
\]
is trivial.
\end{enumerate}

We assume the regularity condition (A.1), since the theory of
global zeta functions is still unsatisfactory for non-regular
prehomogeneous vector spaces (see the final remark (2)).

We denote by $G^+$ the identity component of the real Lie group
$G(\R)$ and put $G^+_x=G^+ \cap G_x$. Let
$\Omega(\R)=\Omega_1\cup\cdots\cup\Omega_\nu$ be the decomposition
into the connected components (in the usual topology). By the
assumption (A.1), every $\Omega_i$ is a single $G^+$-orbit. Let
$\Gamma_N$ be an arithmetic subgroup of $G(\Q)$ which stabilizes
$V(\Z)$ and induces the identity mapping on $V(\Z)/NV(\Z)$. Then
the function $V(\Z) \ni x \mapsto \chi(f(x)) \in \C$ is
$\Gamma_N$-invariant and factors through $V(\Z/n\Z)$. The
$L$-functions $L_i(\mathbf{s};\chi)$ $(1 \leq i \leq \nu)$
associated with $(G,\rho,V)$ and $\chi$ are defined by
\[
L_i(\mathbf{s};\chi) = \sum_{x \in \Gamma_N\backslash (V(\Z)\cap \Omega_i)}
              \mu(x) \chi(f(x))\prod_{j=1}^\ell \abs{P_j(x)}^{-s_j},
              \quad \mathbf{s}=(s_1,\ldots,s_\ell) \in \C^\ell,
\]
where $\mu(x)$ is the volume of the fundamental domain $G^+_x/(\Gamma_N \cap G_x^+)$ with respect to the normalized Haar measure on $G^+_x$ (for the normalization of the Haar measure on $G^+_x$, see \cite{Sa1}, \S4).
By the assumptions (A.1), (A.2) and \cite[Theorem 1.1]{Sai}, the $L$-functions converge absolutely when the real parts of $s_1,\ldots,s_\ell$ are sufficiently large.

We take a relative invariant $f^\vee$ of the dual prehomogeneous
vector space $(G,\rho^\vee,V^\vee)$ with coefficients in $\Z$
(with respect to the basis of $V^\vee$ dual to the fixed basis of
$V$) that corresponds to the character $\phi^{-1}$. Then $f^\vee$
and $(G,\rho^\vee,V^\vee)$ satisfy the assumptions (A.1) and
(A.2). Put $\Omega^\vee = V^\vee \setminus {f^\vee}^{-1}(0)$. We
may order the fundamental relative invariants
$P_1^\vee,\ldots,P_\ell^\vee$ of $(G,\rho^\vee,V^\vee)$ such that
the character corresponding to $P^\vee_i$ is $\phi_i^{-1}$.

Our final assumption is the following:
\begin{enumerate}
\def\labelenumi{(A.\theenumi)}
 \setcounter{enumi}{2}
\item for every prime
factor $p$ of $N$ with $m(p)>1$ (resp.\ $m(p)=1$), Theorems
\ref{theorem A} and \ref{theorem B} in \cite{CH} (resp.\ Theorems
A, B, C in \cite{denef-gyoja}) hold for $\chi^{(p)}$ and $f$.
\end{enumerate}

Let $\psi:\Z/N\Z \rightarrow \C^\times$ be an additive primitive character.
By the Chinese remainder theorem, $\psi$ determines an additive primitive character $\psi^{(p)}:\Z/p^{m(p)}\Z \rightarrow \C^\times$ for each $p|N$.
For $L \in V^\vee(\Z)$, let us consider the following character sum:
\begin{eqnarray*}
S(\chi,f;L)
 &:=& \sum_{x \in V(\Z/N\Z)} \chi(f(x))\psi(L(x)) \\
 &=& \prod_{p|N} \sum_{x \in V(\Z/p^{m(p)}\Z)} \chi^{(p)}(f(x))\psi^{(p)}(L(x)).
\end{eqnarray*}
Then, by the assumptions (A.1) and (A.3), we have
\begin{equation}
\label{eqn:Gauss sum}
S(\chi,f;L)
 = \chi^{-1}(f^\vee(L)) \prod_{p|N} \kappa_p^\vee(L) g(\chi^{(p)},f),
            \end{equation}
with $L \bmod N \in \Omega^\vee(\Z/N\Z)$ and where
$\kappa_p^\vee(L)$ is the constant $\kappa^\vee(L)=\pm1$ given for
each (sufficiently large) $p$ by Theorem \ref{theorem B} in
\cite{CH} or Theorems B and C in \cite{denef-gyoja}, and
$g(\chi^{(p)},f)$ is a constant independent of $L$ whose explicit
value can be easily seen from Theorems \ref{theorem A} and
\ref{theorem B} in \cite{CH} or Theorem A in \cite{denef-gyoja}
according as $m(p)>1$ or $m(p)=1$. Put
\[
\kappa^\vee(L) = \prod_{p|N} \kappa_p^\vee(L) \quad \hbox{and} \quad
g(\chi,f) = \prod_{p|N} g(\chi^{(p)},f).
\]

Now we define the $L$-functions associated with $(G,\rho^\vee,V^\vee)$ by
\[
L_i^\vee(\mathbf{s};\chi^{-1})
 = \sum_{L \in \Gamma_N\backslash (V^\vee(\Z)\cap \Omega^\vee_i)}
              \mu^\vee(L) \kappa^\vee(L)
              \chi^{-1}(f^\vee(L))\prod_{j=1}^\ell \abs{P^\vee_j(L)}^{-s_j}.
\]
The assumptions (A.1), (A.2) and \cite[Theorem 1.1]{Sai} again imply that these $L$-functions converge absolutely when the real parts of $s_1,\ldots,s_\ell$ are sufficiently large.
The abscissa of absolute convergence is independent of $\chi$ and $N$.

To describe analytic properties of the $L$-functions, we need some more notational preliminaries.
Let $b(\mathbf{s})=b(s_1,\ldots,s_\ell)$ be the Bernstein-Sato polynomial defined by
\begin{equation}
\label{eqn:b} \left(\prod_{i=1}^\ell P^\vee_i(\grad_x)\right)
\prod_{i=1}^\ell P_i(x)^{s_i} = b(\mathbf{s})\prod_{i=1}^\ell
P_i(x)^{s_i-1}.
\end{equation}
It is known that the function $b(\mathbf{s})$ is a product of inhomogeneous
linear forms $s_1,\ldots,s_\ell$ of integral coefficients (see \cite{SatoM}).
We also need the Bernstein-Sato polynomial $b_f(s)$ of $f$, which is defined by
\[
f^\vee(\grad_x)f(x)^{s+1} = b_f(s)f(x)^s.
\]
It is known that the roots of $b_f(s)$ are negative rational numbers.

Each of the assumptions (A.1) and (A.2) implies that there exists a
$\delta=(\delta_1,\ldots,\delta_\ell) \in \left(\frac 12\Z\right)^\ell$ such
that the relative invariant $\prod_{i=1}^\ell P_i(x)^{2\delta_i}$ corresponds to
the character $\det\rho(g)^2$ (see \cite[Proposition 8]{SK} or \cite[Proposition
11]{SatoM} for (A.1) and \cite[Lemma 4.1]{Sa1} for (A.2)).

Finally we recall the fundamental theorem of the theory of
prehomogeneous vector spaces over the real number field $\R$. For
$i=1,\ldots,\nu$ and $\mathbf{s}$ with
$\Re(s_1),\ldots,\Re(s_\ell)>0$ , we define a continuous function
$\abs{P(x)}^\mathbf{s}_{\Omega_i}$ on $V(\R)$ by
\[
\abs{P(x)}^\mathbf{s}_{\Omega_i} = \left\{ \begin{tabular}{ll}
                 $\prod_{j=1}^\ell \abs{P_j(x)}^{s_j}$ & $(x \in \Omega_i),$ \cr
                 0 & $(x \not\in \Omega_i).$
                 \end{tabular}\right.
\]
The function $\abs{P(x)}^{\mathbf{s}}_{\Omega_i}$ depends
holomorphically on $\mathbf{s}$ and is extended to a tempered
distribution on $V(\R)$ depending meromorphically on $\mathbf{s}
\in \C^\ell$. We denote the tempered distribution by the same
symbol. We can also define the tempered distributions
$\abs{P^\vee(L)}^\mathbf{s}_{\Omega^\vee_i}$ on $V^\vee(\R)$ quite
similarly. Then the fundamental theorem reads
\begin{equation}
\label{eqn:Fund Thm over R}
\int_{V(\R)} \abs{P(x)}_{\Omega_i}^{\mathbf{s}-\delta} \exp(2\pi i L(x))\,dx
 = \sum_{j=1}^\nu \gamma_{ij}(\mathbf{s}) \abs{P^\vee(L)}_{\Omega^\vee_j}^{-\mathbf{s}},
\end{equation}
where $\gamma_{ij}(\mathbf{s})$ $(i,j=1,\ldots,\nu)$ have
elementary (but not explicit) expressions in terms of the gamma
function and the exponential function (see \cite{MS0},
\cite{Sa1}).

\def\thethm{\Alph{thm}}
\begin{theorem}
\label{thm:L} $(1)$ The $L$-functions $L_i(\mathbf{s};\chi)$ and
$L_i^\vee(\mathbf{s};\chi^{-1})$ multiplied by
$b(\mathbf{s}-\delta)$ have analytic continuations to holomorphic
functions of $\mathbf{s}$ in $\C^\ell$ and satisfy the following
functional equation:
\[
g(\chi,f)
\cdot L^\vee_j(\delta-\mathbf{s};\chi^{-1})
 =  N^{d_1s_1+\cdots+d_\ell s_\ell}\sum_{i=1}^\nu \gamma_{ij}(\mathbf{s}) L_i(\mathbf{s};\chi),
\]
where $d_i$ $(1 \leq i \leq \ell)$ is the degree of the
fundamental relative invariant $P_i$ and $\gamma_{ij}(\mathbf{s})$
is the same as above.

\noindent
$(2)$ Assume that $\chi$ satisfies at least one of the following conditions:
\begin{itemize}
\item $m(p) \geq 2$ for some $p$ dividing the conductor $N$ of
$\chi$; \item  the order of $\chi^{(p)}$ for some $p|N$ with
$m(p)=1$ is different from the reduced denominators of the roots
of the Bernstein-Sato polynomial $b_f(s)$.
\end{itemize}
Then the L-functions $L_i(\mathbf{s};\chi)$ and $L_i^\vee(\mathbf{s};\chi^{-1})$ are entire functions of $\mathbf{s}$ in $\C^\ell$.
\end{theorem}

Since the proof of the theorem is almost the same as the one of
Theorem 2 and Corollary 1 of \cite{Sa1}, we shall give only an
outline of the proof.

Denote by $\BbbA$ the ring of adeles of $\Q$ and by
$\BbbA_0=\prod_{p<\infty}'\Q_p$ the ring of finite adeles of $\Q$.
Denote by $\Phi_p(x_p)$ the characteristic function of $V(\Z_p)$
and put $\Phi_0(x_0)=\prod_{p<\infty} \Phi_p(x_p)$ for
$x_0=(x_p)\in \BbbA_0$. The function $\Phi_0(x_0)$ is the
characteristic function of $\prod_{p<\infty} V(\Z_p)$. Let
$\Phi_\infty(x_\infty)$ is a rapidly decreasing
$C^\infty$-function on $V(\R)$ and define a function $\Phi_{\chi}$
on $V(\BbbA)$ by
\[
\Phi_{\chi}(x) = \Phi_\infty(x_\infty) \cdot \prod_{p|N}
\chi^{(p)}(f(x_p)) \cdot \Phi_0(x_0) \quad (x = (x_\infty,x_0) \in
V(\BbbA)).
\]
The function $\Phi_{\chi}$ is a Schwartz-Bruhat function on
$V(\BbbA)$ and the Poisson summation formula implies the identity
\begin{equation}
\label{eqn:poisson}
\sum_{x \in V(\Q)} \Phi_{\chi}(\rho(g)x)
 = \abs{\det\rho(g)}^{-1}_{\BbbA} \sum_{L \in V^\vee(\Q)}
   \widehat\Phi_{\chi}(\rho^\vee(g)L) \quad (g \in G(\BbbA)),
\end{equation}
where $\widehat\Phi_{\chi}$ is the Fourier transform of
$\Phi_{\chi}$, which is defined by an additive character of
$\BbbA/\Q$ of conductor $1$, more specifically, the additive
character whose $p$-component is of conductor $1$ and coincides
with $\psi^{(p)}(Nx)$ on $p^{-m(p)}\Z_p$ if $p$ divides $N$. It is
easy to see that
\begin{equation}
\label{eqn:FT for fin part}
\widehat\Phi_{\chi}(L) =
N^{-n}S(\chi,f;NL_0)\Phi_0(NL_0)\widehat\Phi_\infty(L_\infty),
\end{equation}
with $L = (L_\infty,L_0) \in V^\vee(\BbbA)$ and $n = \dim V$. We
note here that the function $S(\chi,f;L)$ originally defined on
$V^\vee(\Z)$ can naturally be extended to a function on
$\prod_{p<\infty} V^\vee(\Z_p)$. By the usual technique of
unfolding, we have
\begin{eqnarray}
\lefteqn{
\int_{G^+/\Gamma_N} \prod_{j=1}^\ell \abs{\phi_j(g_\infty)}^{s_j}
 \sum_{x \in \Omega(\Q)} \Phi_{\chi}(\rho(g_\infty,1)x)\,dg_\infty
} \label{eqn:5} \\
  & & = \sum_{i=1}^\nu L_i(\mathbf{s};\chi)
        \langle \abs{P(x)}^{\mathbf{s}-\delta}_{\Omega_i}, \Phi_\infty \rangle. \nonumber
\end{eqnarray}
By the identities (\ref{eqn:Gauss sum}) and (\ref{eqn:FT for fin part}), we also have
\begin{eqnarray*}
\lefteqn{
\int_{G^+/\Gamma_N} \prod_{j=1}^\ell \abs{\phi_j(g_\infty)}^{-s_j}
 \sum_{L \in \Omega^\vee(\Q)}
 \widehat\Phi_{\chi}(\rho^\vee(g_\infty,1)L)\,dg_\infty
}\\
  & & = N^{\sum_{j=1}^\ell d_js_j - n} g(\chi,f)
  \sum_{i=1}^\nu L_i^\vee(\mathbf{s};\chi^{-1})
  \langle \abs{P^\vee(L)}^{\mathbf{s}-\delta}_{\Omega^\vee_i}, \widehat\Phi_\infty \rangle.
\end{eqnarray*}
Now Theorem \ref{thm:L} can be proved in the same manner as in
\cite{Sa1}, \S 6 by using these integral representations of the
$L$-functions, the Poisson summation formula (\ref{eqn:poisson}),
and the fundamental theorem (\ref{eqn:Fund Thm over R}) over $\R$.
We note here that the Poisson summation formula
(\ref{eqn:poisson}) is used in the form
\[
\sum_{x \in \Omega(\Q)} \Phi_{\chi}(\rho(g)x)
 = \abs{\det\rho(g)}^{-1}_{\BbbA} \sum_{L \in \Omega^\vee(\Q)}
   \widehat\Phi_{\chi}(\rho^\vee(g)L) + I(\Phi_\chi, g),
\]
where
\[
I(\Phi_\chi, g)
 = \abs{\det\rho(g)}^{-1}_{\BbbA} \sum_{L \not\in \Omega^\vee(\Q)}
   \widehat\Phi_{\chi}(\rho^\vee(g)L)
   - \sum_{x \not\in \Omega(\Q)} \Phi_{\chi}(\rho(g)x).
\]
The contribution of $I(\Phi_\chi, g)$ to the integral
representation contains the information on the poles of the
$L$-functions and is in general very hard to calculate explicitly.
However, if we take a test function of the form
$\Phi_\infty(x_\infty) = \bigl(\prod_{l=1}^\ell
f^\vee_l(\grad_x)\bigr) \Phi'_\infty(x_\infty)$ for a
$C^\infty$-function $\Phi'_\infty$ with compact support in
$\Omega(\R)$, then $I(\Phi_\chi, g)$ vanishes and, by using
(\ref{eqn:b}), we can prove the first assertion of Theorem
\ref{thm:L}. If the assumption of the second assertion of Theorem
\ref{thm:L} is fulfilled, then Theorem \ref{theorem A} (2) in
\cite{CH} and Remark (5.2.3.3) in \cite{denef-gyoja} imply the
vanishing of $I(\Phi_\chi,(g_\infty,1))$ for $g_\infty \in G^+$
and arbitrary $\Phi_\infty$ with compact support in $\Omega(\R)$.
The second assertion of Theorem \ref{thm:L} follows immediately
from this observation.

\begin{remark}
(1) The simplest example of the $L$-functions considered in this
note is the Dirichlet $L$-function, which is the $L$-function
associated with $(GL_1,\rho,V(1))$, where $\rho$ is the standard 1
dimensional representation of $GL_1$. Further examples of the
$L$-functions associated with prehomogeneous vector spaces have
been studied in \cite{St}, \cite{DW}, \cite{Sai-2}, \cite{Sai-1},
\cite{U1} and \cite{U2}. It is noteworthy that, unlike the case of
the Dirichlet $L$-functions, the $L$-functions $L_i(\bf{s};\chi)$
may have poles even for non-trivial $\chi$, if $\chi$ does not
satisfy any one of the conditions in Theorem \ref{thm:L} (2) (for
concrete examples, see \cite[Theorem 6.2]{DW}, \cite[Theorem
4.2]{Sai-2} and \cite[Theorem 5]{Sai-1}). It is an interesting
problem to determine the conditions on $\chi$ under which
$L_i(\bf{s};\chi)$ has actually poles.

(2) As mentioned before, the theory of global zeta functions is
still incomplete without the assumption (A.1). For example, we do
not have any general convergence theorem for zeta functions. If we
assume the convergence of the zeta integral on the left hand side
of (\ref{eqn:5}), then, we can get a slight generalization of
Theorem \ref{thm:L}. Indeed, by using Theorem 4.19 of
\cite{gyoja}, one can prove an analogue of Theorem \ref{thm:L}
under the following weaker assumption
\begin{enumerate}
\def\labelenumi{(A.\theenumi${}'$)}
 \item Let $O_0$ (resp.\ $O_1$) be the
open (resp.\ closed) orbit in $\Omega$. Then $\Omega=O_0 \cup O_1$
and, for each connected component $\Omega_i$ of $\Omega(\R)$, $O_0
\cap \Omega_i$ is a single $G^+$-orbit.
\end{enumerate}
\end{remark}

\bibliographystyle{plain}

\end{document}